\documentclass[doublespacing]{elsartfrankenstein}

\usepackage{natbib}

 \usepackage{graphicx}
\usepackage{psfrag}
 \usepackage{epsfig,color}

\usepackage{amssymb,amsfonts,amsmath}
\usepackage{latexsym}
\usepackage{mathrsfs}

\newcommand{\be}{\begin{equation}}
\newcommand{\ee}{\end{equation}}
\newcommand{\bea}{\begin{eqnarray}}
\newcommand{\eea}{\end{eqnarray}}
\newcommand{\beas}{\begin{eqnarray*}}
\newcommand{\eeas}{\end{eqnarray*}}
\newcommand{\pa}{\partial}
\newcommand{\parsh}[2]{\frac{\pa #1}{\pa #2}}
\newcommand{\la}{\lambda}
\newcommand{\ve}{\varepsilon}
\newcommand{\ra}{\rightarrow}

\newtheorem{theorem}{Theorem}
\newtheorem{corollary}{Corollary}
\newtheorem{definition}{Definition}

\newcommand{\jerk}{\jmath}
\newcommand{\PPi}{\boldsymbol{\Pi}}

\newcommand{\ihat}{\hat{\mathbf{e}}_x}
\newcommand{\jhat}{\hat{\mathbf{e}}_y}
\newcommand{\khat}{\hat{\mathbf{e}}_z}

\begin{document}

\begin{frontmatter}



\title{The Boltzmann-Hamel Equations for Optimal Control}


\author{Jared M. Maruskin and Anthony M. Bloch}

\address{Department of Mathematics, The University of Michigan}

\begin{abstract}

We will extend the Boltzmann-Hamel equations to the optimal control setting, 
producing a set of equations for both kinematic and dynamic nonholonomic optimal
control problems.  In particular,
we will show the dynamic optimal control problem can be written as a minimal set of $4n-2m$ first
order differential equations of motion.

\end{abstract}


\end{frontmatter}


\section{INTRODUCTION}

\subsection{Overview}  Quasi-velocity formulations, such as Maggi's equation and the Boltzmann-Hamel equation, 
have achieved much success in the analysis of nonholonomic systems due to their ability
to cast the dynamical equations of motion in a form requiring fewer equations, see
 \cite{greenwood2}, \cite{papastavridis}, and \cite{neimark}.
For an $n$ degree of freedom system with $m$ nonholonomic constraints, $2n+m$ equations of motion
are required if one uses the fundamental nonholonomic form of Lagranges equation.  $2n$ differential
equations for the system state, and $m$ algebraic relations that must be solved for the multipliers.  
However, if quasi-velocity techniques are employed, the system can be written as a system of 
$2n-m$ first order differential equations.

The standard approach to optimal control problems is to use Lagrange Multipliers.  Under certain
conditions, the optimal control problem can be reformulated as a vakonomic (variational nonholonomic) problem \cite{blochcrouch}.
One can further analyze optimal control
problems with Pontryagin's Maximum Principle, see \cite{bloch}, \cite{bullo}, or \cite{agrachev}.
Solutions to the kinematic optimal control problems, where one has direct control over a number of the velocities, can be expressed
using $2n+m$ equations of motion; whereas solutions to dynamical optimal control problems, where one has acceleration controls, can be
expressed with $4n+m$ equations of motion.  
Some geometric aspects of this system have been discussed in \cite{camarinha}.
In this paper, we extend quasi-velocity techniques to optimal control problems
with nonholonomic constraints.  We show how to write the optimal control equations for kinematically actuated systems
as a system of $2n$ first order differential equations (a savings of $m$ equations) and the optimal control equations for
dynamically actuated systems as a system of $4n-2m$ first order differential equations (a savings of $3m$ equations).

\subsection{Summation Convention} To aid in notation, we will invoke the 
summation convention throughout this paper.  Greek letters ($\alpha, \beta, \gamma, \ldots$)
run over the constrained dimensions $1, \ldots, m$.  Capital letters ($A, B, C, \ldots$) run
over the unconstrained dimensions $m+1, \ldots, n$.  Lower case letters ($a, b, c, \ldots$) run
over all dimensions $1,\ldots, n$.

\section{QUASI-VELOCITIES AND VARIATIONS}  \label{bhsecaaaII}

In this section we will present the basic background on nonholonomic constraints and
quasi-velocities.    We
will discuss the basic properties of this connection and derive the transpositional
relations, \cite{greenwood2}, \cite{papastavridis}.

\subsection{Nonholonomic Constraints and Quasi-Velocities}

Let $Q$ be the configuration manifold of our system, with $\dim Q = n$ and $TQ$ its corresponding
tangent bundle (our phase space).  A \textit{mechanical} Lagrangian is given by $L : TQ \ra \mathbb{R}$, usually
taken to have the form
$L(q, \dot q) = g_{ij} \dot q^i \dot q^j - V(q) $
where $g_{ij}$ is the \textit{kinetic energy metric} and $V(q)$ is a potential term.

We further suppose our system is subject to $m$ linear scleronomic (time independent) nonholonomic constraints, i.e.
constraints of the form:
\be
\label{bhaaa01}
a^\alpha_i(q) \dot q^i = 0
\ee
Define now a vector space isomorphism $\Psi^j_i$ on the tangent space, with inverse
transformation $\Phi^i_j$.  The first $m$ rows of $\Psi^j_i$ are taken to agree with the
constraint matrix, i.e. $\Psi^\sigma_i(q) = a^\sigma_i(1)$.  The remaining rows can be choosen
freely, so long as the resulting matrix $\Psi$ is invertible.  The transformation
$\Psi$ can be viewed as a change of basis:
\[
\Psi: \left\{ \parsh{}{q^i} \right\}_{i=1}^n \ra \left\{ \parsh{}{\theta^i} \right\}_{i=1}^n
\]
where the new basis is referred to as the \textit{quasi-basis}.  The velocity 
of the system $v \in T_qQ$ can be expressed in terms of the ordinary or quasi-basis
as follows:
\[
\dot q^i \parsh{}{q^i} = \left( \Psi^j_i \dot q^i \right) \parsh{}{\theta^j} = u^j \parsh{}{\theta^j}
= \left( \Phi^i_j u^j \right) \parsh{}{q^i}
\]
where the components $u^j$ are the \textit{quasi-velocities}.  Basis vectors transform as:
\[
\parsh{}{\theta^j} = \Phi^i_j \parsh{}{q^i} \qquad \mbox{and} \qquad \parsh{}{q^i} = 
\Psi^i_j \parsh{}{\theta^j}
\]
Finally, one defines a set of $n$ one-forms, dual to the quasi-basis:
\[
d\theta^j = \Psi^j_i dq^i
\]
Even though this notation is found in the literature, it is really a notational
misnomer, as the one forms $d\theta^j$ are \textit{not} exact.

\subsection{Variations}

\begin{definition} Consider a curve $\gamma(t): [a,b] \ra Q$.  A \textit{proper
variation} of $\gamma(t)$ is a differentiable function $q(s,t):[-\ve, \ve] \times
[a, b] \ra Q$ that satisfies the following conditions:
\begin{itemize}
\item[(i)] $q(0,t) = \gamma(t), \ \ \forall t \in [a,b]$
\item[(ii)] $q(s,a) = \gamma(a)$ and $q(s,b) = \gamma(b), \ \forall s \in [-\ve, \ve]$.
\end{itemize}
\end{definition}

\begin{definition} The \textit{infinitessimal variation} $\delta q(t)$ corresponding to the variation
$q(s,t)$ is the vector field defined along $\gamma(t)$ by
${\displaystyle 
\delta q(t) = \left. \parsh{q(s,t)}{s} \right|_{s=0}}$.
\end{definition}

We will further assume the variations to be \textit{continuous} and \textit{contemporaneous}.  
Continuity of the variations implies that the Lie Derivative $\mathcal{L}_{\dot q} \delta q \equiv 0$ vanishes identically. 
Contemporaneous variations occur without the passage of time.
The infinitessimal variations, when expressed in terms of the quasi-basis, are given by
$ \delta \theta^j(t) = \Psi^j_i \delta q^i$.

\subsection{The Transpositional Relations}

A fundamental ingredient for understanding nonholonomic variational problems
is the following set of transpositional relations (see \cite{greenwood2}, \cite{papastavridis}).

\begin{theorem}[First Transpositional Relations]  Utilizing the shorthand $d := \pa / \pa t$,
$\delta := \pa / \pa s$, we have:
\be
\label{bhaaa02}
(d \delta q^i - \delta d q^i) \Psi^j_i =
(d \delta \theta^j - \delta d \theta^j) + \gamma^j_{ab} u^a \delta \theta^b
\ee
where $\gamma^j_{ab}$ are the Hamel coefficients
${\displaystyle \gamma^s_{pq} = \left\{ \parsh{\Psi^s_i}{q^j} - \parsh{\Psi^s_j}{q^i} \right\} \Phi^i_p \Phi^j_q}$.
\end{theorem}
The left hand side of \eqref{bhaaa02} is no more than $d\theta^j ( \mathcal{L}_{\dot q} \delta q)$;
and, therefore, for continuous variations, is identically zero. We therefore have the following: 

\begin{corollary}  For proper, continuous variations, variations of the quasi-velocities can be related 
to variations of the quasi-coordinates as follows:
\be
\label{bhaaa03}
\delta u^j = d \delta \theta^j + \gamma^j_{ab} u^a \delta \theta^b
\ee
\end{corollary}

Therefore, due to the nonintegrability
of the constraint distribution ($\gamma^\sigma_{ij} \not = 0$, $\sigma = 1, \ldots, m$), one cannot obtain closure in the quasi-coordinate
space, even at the differential
level (\cite{greenwood2}, 
\cite{papastavridis}).  One must choose between $\delta u^\sigma = 0$ or $d \delta \theta^j = 0$.
The correct \textit{dynamical} equations of motion are obtained if one chooses the variations
so that they obey \textit{the Principle of Virtual Work}, $\delta \theta^j \equiv 0$.  
If one, on the other hand, choose the variations to satisfy $\delta u^\sigma=0$, one 
would obtain trajectories that satisfy Hamilton's Principle.  Such trajectories are
referred to as the vakonomic motion of the system, a term introduced by Arnold.

\begin{definition} The associated \textit{quasi-acceleration}, $a^i$, and 
\textit{quasi-jerk}, $\jerk^i$, are defined to be 
$a^i = \dot u^i \qquad \mbox{and} \qquad \jerk^i = \dot{a}^i$.
\end{definition}

A direct coordinate calculation shows:

\begin{theorem}[Second Transpositional Relation] \label{bhaaa04}
For continuous variations,
we have
$\delta d u^i = d \delta u^i$.
Equivalently, $\delta a^i = \pa(\delta u^i)/\pa t$.
\end{theorem}

\section{THE BOLTZMANN-HAMEL EQUATIONS}  \label{bhsecaaaIII}

We will derive the Boltzmann-Hamel equations for nonholonomic mechanics
directly from variational principles.  A more algebraic derivation
of these equations is given in \cite{greenwood2}.  We will begin
with the \textit{Lagrange-D'Alembert Principle}:

\begin{definition}[Lagrange-D'Alembert Principle] The correct \textit{dynamical}
equations of motion are the ones which minimize the action
${\displaystyle I = \int_a^b L(q, \dot q) \ dt}$, 
where $L(q, \dot q)$ is the \textit{unconstrained} mechanical Lagrangian and the variations are chosen to satisfy the Principle of Virtual Work.
\end{definition}

Let $\mathscr{L}(q, u) = L(q, \dot q(q,u))$ be the re-expression of the \textit{unconstrained} Lagrangian
in terms of the quasi-velocities.  Taking variations of the action and using the first transpositional
relations \eqref{bhaaa03}, one obtains:
\beas
\delta I &=& \int_a^b \left( \parsh{\mathscr{L}}{q^i} \delta q^i + \parsh{\mathscr{L}}{u^i} \delta u^i 
+ F_i \delta q^i
\right) \ dt  \\
&=& \int_a^b \left( \parsh{\mathscr{L}}{\theta^i}  - \frac{d}{dt} \parsh{\mathscr{L}}{u^i} 
+ \parsh{\mathscr{L}}{u^j} \gamma^j_{ki} u^k + Q_i \right) \delta \theta^i \ dt
\eeas
where $F_i$ is the external applied force and we have defined:
\[
\parsh{\mathscr{L}}{\theta^i} = \parsh{\mathscr{L}}{q^j} \parsh{q^j}{\theta^i}
= \parsh{\mathscr{L}}{q^j} \Phi^j_i \qquad \mbox{and} \qquad Q_i = \Phi^j_i F_j
\]
After applying the Principle of Virtual Work, $\delta \theta^\sigma \equiv 0$, the remaining
$n-m$ variations $\delta \theta^I$ can be taken to be independent, and we obtain the Boltzmann-Hamel
equations for nonholonomic mechanics:
\bea
\label{bhaaa05}
\frac{d}{dt} \parsh{\mathscr{L}}{u^I} - \parsh{\mathscr{L}}{\theta^I} - \parsh{\mathscr{L}}{u^j}
\gamma^j_{KI} u^K &=& Q_I \\
\label{bhaaa06}
\dot q^i &=& \Phi^i_J u^J
\eea
One must use the \textit{unconstrained} Lagrangian for these equations.  After the partial derivatives
are taken, one then applies the constraints $u^\sigma=0$.  The Boltzmann-Hamel equations \eqref{bhaaa05}-
\eqref{bhaaa06} are a \textit{minimal} set of $2n-m$ first order differential equations
for the $n$ $q^i$'s and the $n-m$ $u^I$'s.

\section{KINEMATIC OPTIMAL CONTROL} \label{bhsecaaaIV}
In this section we will present a quasi-velocity based method for kinematic optimal control 
problems, where one has direct controls over the velocities.  As an example, we will work out
the optimal kinematic control equations for the falling rolling disc.

\subsection{Theory}

For a general affine kinematic control system subject to $m$ nonholonomic constraints, the
following system is typically specified:
$ \dot q^i =  X^i_I(q) w^I$, 
where the $w^I$ are the $n-m$ controls and $X^i_I(q)$ is the $i$-th component of the $I$-th
independent control vector field.  Taking the $m$ constraints as the first $m$ quasi-velocities:
\be
\label{bhaaa08}
u^\sigma = \Psi^\sigma_i \dot q^i \equiv 0,
\ee
one can, wlog, take the controls as the remaining independent quasi-velocities:
\be
\label{bhaaa07}
w^I(q, \dot q) = u^I = \Psi^I_i \dot q^i
\ee
With this choice, the control vector fields are thus identifies with the last
$n-m$ columns of $\Phi = \Psi^{-1}$, i.e.  $X^i_I = \Phi^i_I$.

For a given cost integrand $g(q,w)$, the Kinematic Optimal Control Problem is then given by minimizing the cost function
${\displaystyle
I = \int_a^b g(q, w) \ dt }$
over all curves satisfying \eqref{bhaaa07}-\eqref{bhaaa08} with fixed endpoints $q(a)$ and $q(b)$.

We now define the quasi-basis so that $\Psi^\sigma_i = a^\sigma_i$, as usual, and, additionally,
so that $\Psi^I_i = b^I_i$.  Then the constraints can be written $u^\sigma$, and the $n-m$ 
control variables $w^I$ coincide with the remaining $n-m$ free quasi-velocities $u^I$.  
Define now
$C(q, u) = g(q, w(q, \dot q(q, u)))$.
In our case, we have chosen the unconstrained quasi-velocities to coincide with the controls, i.e. $u^I = w^I$,
thus
we will have $C(q,u) = g(q, u)$.

In order to enforce \eqref{bhaaa08}, 
we must apply the Lagrange Multipliers to the cost function \textit{before} taking variations.  In this
case, we are selecting Hamilton's Principle, where the cost function is minimized amongst the set of
kinematically admissable curves.  We then take unconstrained variations of the augmented cost function
${\displaystyle I = \int_a^b \left( C(q,u) + \mu_\sigma u^\sigma \right) \  dt}$.
Since $C(q,u)$ only depends on the unconstrained quasi-velocities $u^I$, we have:
\[ 
\delta I = \int_a^b \left( \parsh{C}{\theta^i} \delta \theta^i + \parsh{C}{u^I} \delta u^I
+ \mu_\sigma \delta u^\sigma + u^\sigma \delta \mu_\sigma \right) \ dt
\]
Setting the coefficients of the $\delta \mu_\sigma$ terms returns our constraints $u^\sigma = 0$.
Leaving this term off for now, using the transpositional relations \eqref{bhaaa03}, and integrating
by parts yields
\[ 
\delta I = \int_a^b \left\{ \left( \parsh{C}{\theta^i} + \parsh{C}{u^I} \gamma^I_{si} u^s
+ \mu_\sigma \gamma^\sigma_{si} u^s \right) \delta \theta^i 
- \frac{d}{dt} \parsh{C}{u^I} \delta \theta^I
- \dot \mu_\sigma \delta \theta^\sigma \right\} \ dt
\]
We thus have the following 
\begin{theorem} \label{bhtheoremkin} The Boltzmann-Hamel equations for the kinematic optimal control problem are:
\bea
\label{bhaaa09}
\frac{d}{dt} \parsh{C}{u^I} - \parsh{C}{\theta^I} - \parsh{C}{u^J} \gamma^J_{SI} u^S &=& \mu_\tau \gamma^\tau_{SI} u^S \\
\label{bhaaa10}
- \parsh{C}{\theta^\sigma} - \parsh{C}{u^J} \gamma^J_{S\sigma} u^S &=& - \dot \mu_\sigma + \mu_\tau \gamma^\tau_{S\sigma} u^S \\
\label{bhaaa11}
\dot q^i &=& \Phi^i_S u^S
\eea
\end{theorem}
These represent a minimal set of $2n$ first order differential equations: the $n-m$ equations \eqref{bhaaa09} for the unconstrained $u^I$'s,
the $m$ equations \eqref{bhaaa10} for the multipliers $\mu_\sigma$'s, and $n$ kinematic relations \eqref{bhaaa11} for the $q^i$'s.

As an interesting aside, if the cost function integrand $C(q,u)$, when expressed in terms of the 
quasi-velocities, is identical to the \textit{constrained} mechanical Lagrangian, then 
these equations produce the vakonomic motion associated with the system.  See 
\cite{blochcrouch} for additional discussion on the coincidence of the vakonomic motion (Lagrange's Problem)
and the optimal control problem.

\subsection{Optimal Control of the Heisenberg System}

The optimal control of the Heisenberg system, discussed in \cite{brockett} and \cite{bloch},
is a classical underactuated kinematic control problem.  Local coordinates are given by
$q = \langle x, y, z \rangle$.  For this system, one has velocity controls
$w_1 = \dot x$ and 
$w_2 = \dot y$ 
and the motion is subject to the nonholonomic constraint
$\dot z = y \dot x - x \dot y$.
The control velocity field is therefore given by:
\[ \dot q = X_1 w^1 + X_2 w^2, \]
where $X_1 = \langle 1, 0, y \rangle^T$ and $X_2 = \langle 0, 1, -x \rangle^T$.
 Using these controls, one seeks to steer the particle
from the point $\langle 0, 0, 0 \rangle$ at time $t=0$ to the point $\langle 0, 0, a \rangle$
at time $T > 0$, while minimizing the functional
${\displaystyle I = \frac 1 2 \int_0^T \left( w_1^2 + w_2^2 \right) dt}$.

We will derive the equations of motion which yield this solution path via the vakonomic form of the 
Boltzmann-Hamel equations.  We choose quasi-velocities:
$u_1 = y \dot x - x \dot y - \dot z$, $\ u_2 = \dot x$, and $u_3 = \dot y$.  
Notice the quasi-velocities $u_2$ and $u_3$ coincide with the control velocities.
The transformation matrices $\Psi$ and $\Phi$ are given by: 
\[
\Psi = \left[ \begin{array}{ccc}
 y & -x & -1 \\
1 & 0 & 0 \\
0 & 1 & 0 \end{array} \right]
\ \ \mbox{and} \ \
\Phi = \left[ \begin{array}{ccc}
0 & 1 & 0 \\
0 & 0 & 1 \\
-1 & y & -x \end{array} \right]
\]
The nonzero Hamel coefficients are
$ \gamma^1_{23} = - \gamma^1_{32} = 2$.  
Expressing the integrand of the cost function in terms of quasi-velocities yields
${\displaystyle C = \frac 1 2 \left( u_2^2 + u_3^2 \right)}$.  
The kinematic optimal control Boltzmann-Hamel equations \eqref{bhaaa09}-\eqref{bhaaa11} immediately produce the following
set of first order differential equations:
\beas
\dot x = u_2 \qquad & \dot y = u_3 & \qquad \dot z = -u_2 + yu_2 - x u_3 \\
\dot u_2 = -2 \mu u_3 \qquad & \dot u_3 = 2 \mu u_2 & \qquad \dot \mu = 0
\eeas
where $\mu(t) = \mu(0)$ is an arbitrary constant that can be choosen such that 
the solution trajectory reaches its final destination point.  The top equations are a reiteration
of the control field $\dot q = X_1 w^1 + X_2 w^2 = X_1 u^2 + X_2 u^3$ and the bottom equations
produce the optimal control.

\subsection{Optimal Control of the Vertical Rolling Disc} \label{bhsecVRD}

The generalized coordinates of the vertical rolling disc are given by $q = \langle x, y, \theta, \phi \rangle$,
where $(x, y)$ is the contact point of the disc and the $x-y$ plane, $\phi$ is the angle the disc makes with
the $x$-axis, and $\phi$ is the angle a reference point on the disc makes with the vertical.  
Assume we have the kinematic controls
$w_1 = \dot \theta$ and $w_2 = \dot \phi$,  
and that the motion is subject to the nonholonomic constraints
$\dot x - \cos(\phi) \dot \theta = 0$ and 
$\dot y - \sin(\phi) \dot \theta = 0$. 
This gives rise to the control vector field
$\dot q = X_1 w^1 + X_2 w^2$ 
where $X_1 = \langle \cos \phi, \sin \phi, 1, 0 \rangle^T$ and $X_2 = \langle 0, 0, 0, 1 \rangle^T$.

We wish to steer the disc between two points while minimizing the cost functional
${\displaystyle \frac 1 2 \int_a^b (w_1^2 + w_2^2) \ dt}$.  
We choose quasi-velocities
$u_1 = \dot x - \cos(\phi) \dot \theta$, 	
$u_2 = \dot y - \sin(\phi) \dot \theta$, $u_3 = \dot \theta$, and $u_4 = \dot \phi$, 
so that the transformation matrices $\Psi$ and $\Phi$ are given by:
\[
\Psi = \left[ \begin{array}{cccc}
1 & 0 & -\cos \phi & 0 \\
0 & 1 & -\sin \phi & 0 \\
0 & 0 & 1 & 0 \\
0 & 0 & 0 & 1
\end{array} \right] \ \
\Phi = \left[ \begin{array}{cccc}
1 & 0 & \cos \phi & 0 \\
0 & 1 & \sin \phi & 0 \\
0 & 0 & 1 & 0 \\
0 & 0 & 0 & 1
\end{array} \right]
\]

The Hamel coefficients are:
$\gamma^1_{34} = \sin \phi = - \gamma^1_{43}$ and $\gamma^2_{34} = -\cos \phi = -\gamma^2_{43}$.
In terms of the quasi-velocities, the integrand of the cost function  becomes
$C(q,u) = \frac 1 2 u_3^2 + \frac 1 2 u_4^2$. 
The  Boltzmann-Hamel equations \eqref{bhaaa09}-\eqref{bhaaa11} then produce the following set of first order differential equations:
\beas
\dot u_3 = (\mu_2 \cos \phi - \mu_1 \sin \phi) u_4 \qquad \dot \mu_1 = 0 & \qquad \dot x = \cos(\phi) u_3 \qquad & \dot \theta = u_3 \\
\dot u_4 = (\mu_1 \sin \phi - \mu_2 \cos \phi) u_3 \qquad \dot \mu_2 = 0 & \dot y = \sin(\phi) u_3 & \dot \phi = u_4
\eeas

\subsection{Kinematic Optimal Control of the Falling Rolling Disc}

\begin{figure}[h]
\begin{center}
\psfrag{ed}{$\e_d$}
\psfrag{dphi}{$\dot{\boldsymbol{\phi}}$}
\psfrag{dtheta}{$\dot{\boldsymbol{\theta}}$}
\psfrag{dpsi}{$\dot{\boldsymbol{\psi}}$}
\psfrag{etheta}{$\e_\theta$}
\psfrag{epsi}{$\e_\psi$}
\psfrag{P}{$P$}
\psfrag{C}{$C$}
\psfrag{r}{$r$}
\psfrag{x}{$x$} \psfrag{y}{$y$}
\psfrag{z}{$z$} \psfrag{phi}{$\phi$}
\psfrag{psi}{$\psi$}
\psfrag{theta}{$\theta$}
\includegraphics[width=2.7in]{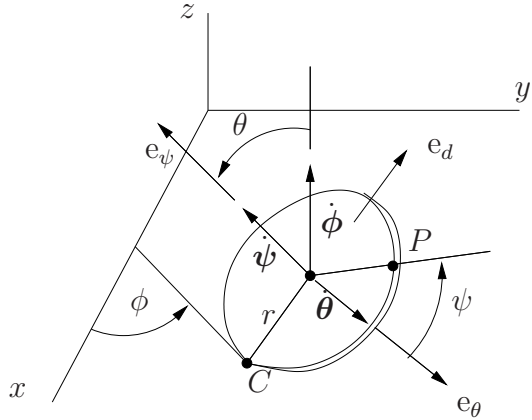}
\end{center}
\caption{Euler Angles of the Falling Rolling Disc}
\label{bhoptcon/figa06}
\end{figure}

The falling rolling disc can be described by the contact point $(x,y)$ 
and Classical Euler angles $(\phi, \theta, \psi)$, as shown in Figure \ref{bhoptcon/figa06}.
We will take the coordinate ordering $(\phi, \theta, \psi, x, y)$.

Suppose we have direct control over the body-axis angular velocities
$w_1 = \omega_d := \dot \phi \sin \theta$, $w_2 = \dot \theta$, and  $w_3 = \Omega := \dot \phi \cos \theta + \dot \psi$
(in the $\e_d$, $\e_\theta$, and $\e_\psi$ directions, respectively (see Fig. \ref{bhoptcon/figa06})),
and the system is subject to the nonholonomic constraints
$\dot x + r \dot \psi \cos \phi = 0$ and $\dot y + r \dot \psi \sin \phi = 0$.
We wish to steer the disc between two points while minimizing the cost functional
${\displaystyle I[\gamma] = \frac 1 2 \int_a^b \left( w_1^2 + w_2^2 + w_3^2 \right) \ dt}$.
We will choose as quasi-velocities
$u_1 = \dot \phi \sin \theta$, $u_2 = \dot \theta$,  $u_3 = \dot \phi \cos \theta + \dot \psi$, 
$u_4 = \dot x + r \dot \psi \cos \phi$, and $u_5 = \dot y + r \dot \psi \sin \phi$.
The quasi-velocities $(u_1, u_2, u_3) = (\omega_d, \dot \theta, \Omega)$
represent the angular velocity expressed in the body-fixed frame, and are coincident with 
the kinematic controls.  These 
are not true velocities (like the Euler Angle Rates), as they are non-integrable.
The nonholonomic constraints in terms of these variables are $u_4 = u_5 = 0$.

The transformation matrices
are
\[
\Psi = \left[ \begin{array}{ccccc}
\sin \theta & 0 & 0 & 0 & 0 \\
0 & 1 & 0 & 0 & 0 \\
\cos \theta & 0 & 1 & 0 & 0 \\
0 & 0 & r \cos \phi & 1 & 0 \\
0 & 0 & r \sin \phi & 0 & 1
\end{array} \right]
\qquad \mbox{and} \qquad
\Phi = \left[ \begin{array}{ccccc}
\csc \theta & 0 & 0 & 0 & 0 \\
0 & 1 & 0 & 0 & 0 \\
-\cot \theta & 0 & 1 & 0 & 0 \\
r \cos \phi \cot \theta & 0 & -r \cos \phi & 1 & 0 \\
r \sin \phi \cot \theta & 0 & -r \sin \phi & 0 & 1
\end{array} \right]
\]

The nonzero Hamel-coefficients are
$\gamma^1_{21} = -\cot \theta = - \gamma^1_{12}$, $\gamma^3_{21} =  1 = - \gamma^3_{12}$,
$\gamma^4_{13} =  r \sin \phi \csc \theta = - \gamma^4_{31}$, and
$\gamma^5_{13} = -r \cos \phi \csc \theta = - \gamma^5_{31}$.

Written in terms of the quasi-velocities, the integrand of the cost function is 
$C(q, u) = \frac 1 2 (u_1^2 + u_2^2 + u_3^2)$.
The kinematic optimal control Boltzmann-Hamel equations \eqref{bhaaa09}-\eqref{bhaaa11} give us a minimal set
of 10 first order differential equations:
\beas
\dot u_1 &=& u_2 u_3 -  u_1 u_2 \cot \theta - r( \mu_4 \sin \phi - \mu_5 \cos \phi) \csc \theta u_3 \\
\dot u_2 &=& u_1^2 \cot \theta - u_1 u_3 \\
\dot u_3 &=& r (\mu_4 \sin \phi - \mu_5 \cos \phi) \csc \theta u_1 \\
\dot \mu_4 &=& 0, \ \ 
\dot \mu_5 = 0 \\
\dot \phi &=& \csc \theta u_1, \ \
\dot \theta = u_2, \ \ 
\dot \psi = - \cot \theta u_1 + u_3 \\
\dot x &=& r \cos \phi \cos \theta u_1 - r \cos \phi u_3, \ \ 
\dot y = r \sin \phi \cot \theta u_1 - r \sin \phi u_3 
\eeas

\section{DYNAMIC OPTIMAL CONTROL}  \label{bhsecaaaV}
In this section, we will derive a set of Boltzmann-Hamel equations for the dynamic optimal control
problem, which is normally a fourth order system.  We will present a minimal set of $4n-2m$ first
order differential equations that produces the optimal control, and then discuss examples.

\subsection{Boltzmann-Hamel Equations for Optimal Dynamic Control}

Given a nonholonomic mechanical system with $n-m$ independent acceleration controls, it
can be recast into the form given by the dynamical Boltzmann-Hamel equations \eqref{bhaaa05}-\eqref{bhaaa06}.
The dynamical optimal control problem is the problem of finding solution curves between
two fixed points $\langle q(a), \dot q(a) \rangle$ and $\langle q(b), \dot q(b) \rangle$ that minimize the cost
function
${\displaystyle I = \int_a^b g(q, \dot q, Q) \ dt}$.  
Utilizing \eqref{bhaaa05} and \eqref{bhaaa06}, we can rewrite the integrand as an explicit function of the coordinates,
quasi-velocities, and quasi-accelerations
$C(q, u, a) = g(q, \dot q(q, u), Q(q, u, a))$.

Since the Boltzmann-Hamel equations no longer depend on the constrained quasi-velocities and quasi-accelerations,
$C(q,u, a)$ is also independent of $u^\sigma$ and $a^\sigma$.  Taking variations
yields:
${\displaystyle \delta I = 
\int \left\{ \parsh{{C}}{q^i} \delta q^i + \parsh{{C}}{u^J} \delta u^J
+ \parsh{{C}}{a^J} \delta a^J \right\} dt}$.  
Using the second transpositional relations Theorem \ref{bhaaa04}
for $\delta a^J$
and then integrating by parts
we obtain
${\displaystyle 
\delta I =
\int \left\{ \parsh{{C}}{q^i} \delta q^i + \left[ \parsh{{C}}{u^J} - \frac{d}{dt} \parsh{{C}}{a^J}
\right] \delta u^J \right\} dt}$.
Defining the parameters
\be
\label{bhaaa12}
\kappa_J = \parsh{{C}}{u^J} - \frac{d}{dt} \parsh{{C}}{a^J}
\ee
and using the first Transpositional relations
\eqref{bhaaa03}  we obtain:
\[
\delta I = \int \left\{ \parsh{{C}}{\theta^r} - \dot \kappa_J \delta^J_r 
+ \kappa_J \gamma^J_{sr} u^s \right\} \delta \theta^r
\ dt
\]
These variations are not free, but subject to the nonholonomic constraints $a^\sigma_i \dot q^i = 0$.  We form the augmented
cost integrand by replacing $C(q,u,a)$ with $C(q, u, a) + \mu_\sigma u^\sigma$.  Taking variations,
the $\delta \mu^\sigma$ coefficients recover the constraints.  Ignoring these terms, we are left with
$\delta I =$
\[ 
 \int \left\{ \parsh{{C}}{\theta^r} - \dot \kappa_J \delta^J_r 
+ \kappa_J \gamma^J_{sr} u^s 
 - \dot \mu_\sigma \delta^\sigma_r + \mu_\sigma \gamma^\sigma_{sr} u^s \right\} \delta \theta^r \ dt
\]
where the variations are now taken to be unconstrained.  Notice the multipliers $\mu_\sigma$ are \textit{not}
the \textit{mechanical} multipliers, but a multiplier on the cost function that enforces Hamilton's Principle.
We thus have the following:
\begin{theorem} \label{bhtheoremdyn} The Boltzmann-Hamel equations for Optimal Dynamic Control are given by:
\bea
\label{bhaaa13}
- \parsh{{C}}{\theta^A} + \dot{\kappa}_A - \kappa_J \gamma^J_{SA} u^S &=&  \mu_\tau \gamma^\tau_{SA} u^S \\
\label{bhaaa14}
- \parsh{{C}}{\theta^\sigma} - \kappa_J \gamma^J_{S \sigma} u^S &=&  \mu_\tau \gamma^\tau_{S \sigma} u^S - \dot \mu_\sigma \\
\label{bhaaa24}
\dot q^i &=& \Phi^i_S u^S
\eea
\end{theorem}
The optimal control system can therefore be given by a minimal set of $4n-2m$ first order differential equations as follows.
We have $n$ kinematic relations \eqref{bhaaa24}, $2n-2m$ relations
$\dot u^A = a^A$ and $\dot a^A = \jerk^A$, 
$n-m$ equations for  $\dot{\jerk}^A$ (given by inserting \eqref{bhaaa12} into \eqref{bhaaa13}), and, finally, $m$ relations
for the multipliers $\dot{\mu}_\sigma$ \eqref{bhaaa14}.   Once the resulting optimal control dynamics
are determined, the control forces which produce the optimal trajectory are then given by the $n-m$ algebraic equations
\eqref{bhaaa05}.  The solution is then found by solving the related boundary value problem, with $4n-2m$ prescribed
boundary conditions:  $q^i(0), \ u^A(0), \ q^i(T), \ u^A(T)$.

\subsection{Dynamic Optimal Control of the Vertical Rolling Disc}

Consider the vertical rolling disc of \S \ref{bhsecVRD} with control torques in the $\theta$ and $\phi$
directions. 
The corresponding dynamical equations of motion (see \cite{bloch}) are:
$\frac 3 2 \ddot \theta = w_3$,  $\frac 1 4 \ddot \phi = w_4$, $\dot x = \dot \theta \cos \phi$, and $\dot y = \dot \theta \sin \phi$.
This is equivalent to a minimal set of 6 first order differential equations (the number
obtained by using the Boltzmann-Hamel equations \eqref{bhaaa05} and \eqref{bhaaa06}.

We now wish to choose the control forces so as to minimize the cost function
$\int \frac 1 2 (w_3^2 + w_4^2) \ dt$.
Solving for the controls in terms of the quasi-accelerations $w_3 = \frac 3 2 \ddot \theta =
\frac 3 2 a_3$ and $w_4 = \frac 1 4 \ddot \phi = \frac 1 4 \ddot a_4$, this is equivalent to
minimizing the action
${\displaystyle \int \left( \frac 9 8 a_3^2 + \frac 1 {32} a_4^2 \right) \ dt}$
subject to the nonholonomic constraints.  Using the dynamic optimal control Boltzmann-Hamel equations
\eqref{bhaaa13} and \eqref{bhaaa14}, coupled with the dynamical equations of motion above, and eliminating the 
controls, we have a minimal system of 12 first order differential equations:
\beas
\dot x = \cos \phi \ u_3 & \qquad & \dot{\jerk}_3 = \frac 4 9 (\mu_1 \sin \phi - \mu_2 \cos \phi) u_4 \\
\dot y = \sin \phi \ u_3 & & \dot{\jerk}_4 = 16(- \mu_1 \sin\phi + \mu_2 \cos \phi) u_3 \\
\dot \theta = u_3 & \dot u_3 = a_3 & \dot a_3 = \jerk_3 \ \ \ \  \ \dot \mu_1 = 0 \\
\dot \phi = u_4  & \dot u_4 = a_4 & \dot a_4 = \jerk_4 \ \ \ \ \ \dot \mu_2 = 0
\eeas
By use of quasi-velocities, quasi-accelerations,
and quasi-jerks, we have made the following simplifications: $u_1 = u_2 = a_1 = a_2 = \jerk_1 = \jerk_2 = 0$,
thereby eliminating the necessity of 6 of the 18 first order differential equations necessary in the standard
approach.
The solution to this system of differential equations yields the optimal dynamic control equations 
of the vertical rolling disc.  It is equivalent to the following reduced system
\beas
\dot x = \cos \phi \dot \theta & \qquad & \dot y = \sin \phi \dot \theta \\
\ddddot \theta = \frac 4 9 (\mu_1 \sin \phi - \mu_2 \cos \phi) \dot \phi & & 
\ddddot \phi = 16(- \mu_1 \sin \phi + \mu_2 \cos \phi) \dot \theta
\eeas
where $\mu_1, \ \mu_2$ are constants.

\subsection{Dynamic Optimal Control of the Free Rigid Body}

Consider dynamic control of the free rigid body, where the generalized coordinates are given by the Type-I Euler
angles $(\psi, \theta, \phi)$.  As quasi-velocities, choose the body-fixed components of the angular momentum
$u_1 = \omega_x = - \dot \psi \sin \theta + \dot \phi$, 
$u_2 = \omega_y = \dot \psi \cos \theta \sin \phi + \dot \theta \cos \phi$, and 
$u_3 = \omega_z = \dot \psi \cos \theta \cos \phi - \dot \theta \sin \phi$.
The transformation matrices are given as:
\[
\Psi = \left[ \begin{array}{ccc}
- \sin \theta & 0 & 1 \\
\cos \theta \sin \phi & \cos \phi & 0 \\
\cos \theta \cos \phi & - \sin \phi & 0
\end{array} \right]
\qquad \mbox{and} \qquad
\Phi = \left[ \begin{array}{ccc}
0 & \sec \theta \sin \phi & \sec \theta \cos \phi \\
0 & \cos \phi & - \sin \phi \\
1 & \tan \theta \sin \phi & \tan \theta \cos \phi
\end{array} \right]
\]
The mechanical Lagrangian is given as
${\displaystyle 
\mathscr{L}(q,u) = \frac 1 2 (I_{xx} u_1^2 + I_{yy} u_2^2 + I_{zz} u_3^2 )}$.  
The nonzero Hamel coefficients are
$\gamma^1_{23} = 1$, $\gamma^2_{13} =  -1$, $\gamma^3_{12} = 1$, 
$\gamma^1_{32} = - 1$, $\gamma^2_{31} = 1$, and $\gamma^3_{21} = -1$.
For notational convenience, define
$\eta_{32} = I_{zz} - I_{yy}$, 
$\eta_{13} = I_{xx} - I_{zz}$, and
$\eta_{21} = I_{yy} - I_{xx}$.
Then the Boltzmann-Hamel equations \eqref{bhaaa05} produce the Euler Equations:
\be
\label{bhaaa17}
I_{xx} \dot u_1 + \eta_{32} u_2 u_3 = M_x \qquad
I_{yy} \dot u_2 + \eta_{13} u_1 u_3 = M_y \qquad
I_{zz} \dot u_3 + \eta_{21} u_1 u_2 = M_z
\ee
where $M_x$, $M_y$, and $M_z$ are the control torques applied about the body fixed principal axes.
The cost function integrand $\frac 1 2 (M_x^2 + M_y^2 + M_z^2)$, when expressed in terms of quasi-variables, is given by:
$C =  \frac 1 2 \{ I_{xx}^2 a_1^2 + I_{yy}a_2^2 + I_{zz}a_3^2 + 2 I_{xx} \eta_{32} a_1 u_2 u_3 
+ 2 I_{yy} \eta_{13} u_1 a_2 u_3 + 2 I_{zz}\eta_{21}u_1 u_2 a_3 
+ \eta_{32}^2 u_2^2 u_3^2 + \eta_{13}^2 u_1^2 u_3^2 + \eta_{21}^2 u_1^2 u_2^2 \}$.

The $\kappa$'s \eqref{bhaaa12} are given by:
\bea
\label{bhaaa27}
\kappa_1 &=& I_{yy} \eta_{13} a_2 u_3 +  I_{zz} \eta_{21} u_2 a_3 
 +  \eta_{13}^2 u_1 u_3^2 +  \eta_{21}^2 u_1 u_2^2  \\
& & -  I_{xx} \jerk_1 -  I_{xx} \eta_{32} u_2 a_3 -  I_{xx} \eta_{32} a_2 u_3 \nonumber \\
\kappa_2 &=& I_{xx} \eta_{32} a_1 u_3 +  I_{zz} \eta_{21} u_1 a_3 
 +  \eta_{32}^2 u_2 u_3^2 +  \eta_{21}^2 u_1^2 u_2  \\
& & -  I_{yy} \jerk_2 -  \eta_{13} I_{yy} u_1 a_3 -  \eta_{13} I_{yy} a_1 u_3 \nonumber \\
\label{bhaaa29}
\kappa_3 &=& I_{xx} \eta_{32} a_1 u_2 + I_{yy} \eta_{13} u_1 a_2 
 +  \eta_{32}^2 u_2^2 u_3 +  \eta_{13}^2 u_1^2 u_3 \\
& & -  I_{zz} \jerk_3 -  \eta_{21} I_{zz} u_1 a_2 -  \eta_{21} I_{zz} a_1 u_2 \nonumber
\eea
The optimal control Boltzmann-Hamel
equations \eqref{bhaaa13} then work out to be:
\be
\label{bhaaa20}
\dot{\boldsymbol{\kappa}} = \boldsymbol{\kappa} \times \boldsymbol{\omega}
\ee
These provide 3 differential equations for the $\dot{\jerk}$'s.  Let $\mathbb{I}$ be the 
moment inertia tensor with respect to the principal axes basis $\ihat, \ \jhat, \ \khat$,
so that, in dyadic notation, $\mathbb{I} = I_{xx} \ihat \ihat + I_{yy} \jhat \jhat
+ I_{zz} \khat \khat$.  Let  $\PPi := \mathbb{I} \cdot \boldsymbol{\omega}$
be the body axis angular momentum, and $\boldsymbol{\kappa} = \langle
\kappa_1, \kappa_2, \kappa_3 \rangle$.  Then \eqref{bhaaa27}-\eqref{bhaaa29} can alternatively be re-expressed
as:
\be
\label{bhaaa30}
\boldsymbol{\kappa} = \PPi \times \dot{\PPi} + \PPi \times (\boldsymbol{\omega} \times \PPi)
- \ddot{\PPi} 
 - \mathbb{I} \cdot \left\{
2 \boldsymbol{\omega} \times \dot{\PPi} + \dot{\boldsymbol{\omega}} \times \PPi +
\boldsymbol{\omega} \times (\boldsymbol{\omega} \times \PPi) \right\}  
\ee
Finally, by defining $\boldsymbol{\la}(\boldsymbol{\omega}, \dot{\boldsymbol{\omega}}) = \boldsymbol{\kappa} + \ddot{\PPi}$, the dynamic optimal control equations
for the free rigid body can be expressed as:
\be
\label{bhaaa31}
\dddot{\PPi} = \dot{\boldsymbol{\la}} + \ddot{\PPi} \times \boldsymbol{\omega} - 
\boldsymbol{\la} \times \boldsymbol{\omega}
\ee
In addition, we have the kinematic relations
\bea
\label{bhaaa23}
\dot \psi &=& \sec \theta \sin \phi u_2 + \sec \theta \cos \phi u_3 \\
\dot \theta &=& \cos \phi u_2 - \sin \phi u_3 \\
\label{bhaaa25}
\dot \phi &=& u_1 + \tan \theta \sin \phi u_2 + \tan \theta \cos \phi u_3
\eea
as well as the relations $\dot u_i = a_i, \ \dot a_i = \jerk_i$.  This is a set of 12 first order 
differential equations.  Once one solves the corresponding boundary value problem
(initial, final Euler angles, angular velocities specified), the controls are determins by the algebraic relations \eqref{bhaaa17}.

\begin{figure}[t]
\begin{center}
\psfrag{A}{$\psi$}
\psfrag{B}{$\theta$}
\psfrag{C}{$\phi$}
\psfrag{E}{$\omega_x$}
\psfrag{F}{$\omega_y$}
\psfrag{G}{$\omega_z$}
\includegraphics[width=3in]{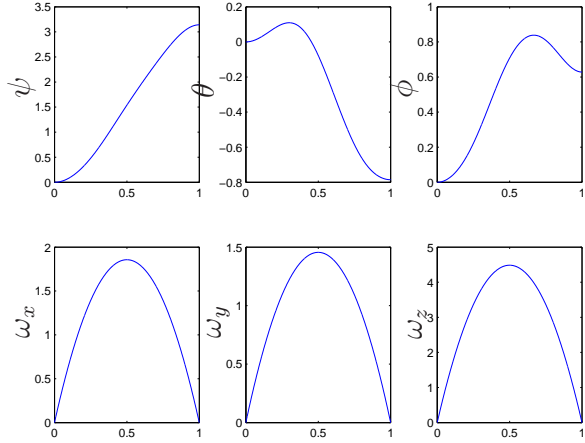}
\end{center}
\caption{Optimal Dynamic Control of Free Sphere:  Euler Angles and Body Fixed Angular Velocity with respect to time.}
\label{bhoptcon/figaaa01}
\end{figure}

For the special case when the rigid body is spherical one sees
from \eqref{bhaaa30} that $\boldsymbol{\kappa} = - \ddot{\boldsymbol{\Pi}}$ and $\boldsymbol{\la} = \mathbf{0}$.  Then the Boltzmann-Hamel 
equations for the optimal dynamic control of the free rigid body \eqref{bhaaa31} reduce to
$\dddot{\boldsymbol{\omega}} = \ddot{\boldsymbol{\omega}} \times {\boldsymbol{\omega}}$.  
When coupled with the kinematic relations \eqref{bhaaa23}-\eqref{bhaaa25} and the algebraic relations \eqref{bhaaa17},
the optimal control trajectories of the free rigid sphere are produced.  
Integrating once yields the second order system
$\ddot{\boldsymbol{\omega}} = \mathbf{c} + \dot{\boldsymbol{\omega}} \times {\boldsymbol{\omega}}$, 
which coincides with the result of \cite{noakes}.  See also \cite{crouchleite}.
The optimal solution trajectory of the reorientation of the rigid sphere from $\mathbf{q}(0) = \langle 0, 0, 0 \rangle,
\ \boldsymbol{\omega}(0) = \langle 0, 0, 0 \rangle$ to the point $\mathbf{q}(1) = \langle \pi, -\pi/4, \pi/5 \rangle, \
\boldsymbol{\omega}(1) = \langle 0, 0, 0 \rangle$ is plotted in Fig. \ref{bhoptcon/figaaa01}.

\section{CONCLUSIONS}

In this paper, we showed how one can extend quasi-velocity techniques to kinematic and optimal
control problems.  Standard Lagrange Multiplier techniques for kinematical optimal control problems
produce a set of $2n+m$ first order differential equations:  $n$ for the coordinates $q^i$, $n$ for
the velocities $\dot q^i$, and $m$ for the multipliers $\mu_\sigma$.  On the other hand, by generalizing the dynamic Boltzmann-Hamel
equations to the kinematic control setting (Theorem \ref{bhtheoremkin}), we obtain a savings of $m$ first order differential equations,
as one no longer need solve for the constrained quasi-velocities. 
Moreover, the differential equations for the multipliers \eqref{bhaaa10} are naturally separated from the differential
equations for the quasi-velocities \eqref{bhaaa09}.

For the dynamic optimal control problem, one typically encounters a fourth order system, plus multipliers, which produces
a total of $4n+m$ first order differential equations.  The Boltzmann-Hamel form of the equations (Theorem \ref{bhtheoremdyn})
gives a minimal set of $4n-2m$ equations of motion, 
as one no longer need integrate the $m$ constrained quasi-velocities, quasi-accelerations, 
and quasi-jerks, $u^\sigma \equiv 0$, $a^\sigma \equiv 0$, $\jerk^\sigma \equiv 0$, respectively.  This approach
gives us a total savings of $3m$ first order differential equations.
Initial and final conditions are then enforced
by solving the resulting system of differential equations as a two point boundary value problem.

The authors wish to thank support from NSF grants DMS-0604307 and CMS-0408542.




\end{document}